\documentclass{amsart}
\usepackage{amscd, amssymb, amsmath, amsthm}
\usepackage{amsfonts,enumerate}
\usepackage{graphics,graphicx,psfrag}
\usepackage[colorlinks=true]{hyperref}

\newtheorem{theorem}{Theorem}[section]
\newtheorem{lemma}[theorem]{Lemma}
\newtheorem{prop}[theorem]{Proposition}

\theoremstyle{definition}
\newtheorem{definition}[theorem]{Definition}
\newtheorem{convention}[theorem]{Convention}

\newtheorem{remark}[theorem]{Remark}

\newcommand{\QQ}{{\mathbb Q}}
\newcommand{\CC}{{\mathbb C}}
\newcommand{\ZZ}{{\mathbb Z}}
\newcommand{\FF}{{\mathbb F}}
\newcommand{\HH}{{\mathbb H}}

\newcommand{\seq}[1]{{\{ {#1}\}}}

\newcommand{\SL}{{\mathrm{SL}}}
\newcommand{\mobgp}{{\mathrm{PSL}_2(\mathbb{C})}}
\newcommand{\Isom}{{\mathrm{Isom}}}
\newcommand{\Hom}{{\mathrm{Hom}}}
\newcommand{\id}{{\mathrm{id}}}

\title[Homomorphism J{\o}rgensen-Thurston]{A J{\o}rgensen-Thurston theorem for homomorphisms}
\date{\today}

\author[Yi Liu]{Yi Liu} 
\address{
    University of California, Berkeley \\
    970 Evans Hall \#3840 \\
    Berkeley, CA 94720-3840} 
\email{yliu@math.berkeley.edu}

\begin{document}


\begin{abstract}
In this note, we provide a description of the structure
homomorphisms from a finitely generated group to 
any torsion-free ($3$-dimensional) Kleinian group with uniformly bounded
finite covolume. This is analogous to the J{\o}rgensen-Thurston Theorem
in hyperbolic geometry.
\end{abstract}

\maketitle

\section{Introduction}\label{Sec-Intro}
	In classical $3$-dimensional hyperbolic geometry, the J{\o}rgensen-Thurston Theorem (\cite[Theorem
	5.12.1]{Th}) asserts that for every $V>0$,
	every complete orientable hyperbolic $3$-manifold of volume at most $V$
	is isometric to a hyperbolic Dehn filling of one of 
	a finite collection of complete
	orientable hyperbolic $3$-manifolds of volume at most $V$. 
	This can be restated using Kleinian groups (cf.~Theorem \ref{JT}). We always speak
	of Kleinian groups in the classical sense, meaning
	discrete subgroups of $\mobgp$. 

	In this note, we prove the following theorem:
	
	\begin{theorem}\label{main-JTh}
		Suppose $G$ is a finitely generated group, and 
		suppose $V>0$. Then there exist finitely many groups $R_1,\cdots, R_k$,  
		each of which is a Dehn extension of a torsion-free Kleinian group of covolume
		at most $V$, such that for any torsion-free Kleinian group $H$ of covolume at most $V$, 
		there is some $1\leq i\leq k$, so that every homomorphism $\phi:G\to H$
		factors through $R_i$. Namely, there is some $\psi:G\to R_i$ such that
		$\phi=\iota_{i,H}\circ\psi$, where $\iota_{i,H}:R_i\twoheadrightarrow H$
		is an extended-filling epimorphism associated to $R_i$ and $H$, independent from $\phi$.
	\end{theorem}
		
	As $R_i$'s are all torsion-free hyperbolic groups relative to isolated abelian subgroups,
	(see the explanation after Proposition \ref{extStable}),
	the structure of homomoprhisms in $\Hom(G,R_i)$ can be described using the Makanin-Razborov diagram
	by the work of Daniel Groves (\cite{Gr-HypII}, cf.~\cite{Se} for an original version
	for free groups). In this sense, Theorem \ref{main-JTh} implies a uniform
	description of the structure of homomorphisms from a finitely generated group to torsion-free Kleinian groups
	of uniformly-bounded volume.
	
	When $G$ is finitely presented, Theorem \ref{main-JTh} is a quick consequence of the factorization
	result obtained in \cite[Theorem 3.2]{AL}, (Section \ref{Sec-fpCase}). 
	To obtain the finitely generated case, we show that 
	generic Dehn extensions	of finite-covolume torsion-free Kleinian groups 
	are stable limits of congruent extended fillings, (Proposition \ref{extStable}).
	This is also a special case of some more general results 
	for Dehn fillings of relatively hyperbolic groups, (cf.~Subsection \ref{Subsec-extStable}).
	It follows from a trick of representation varieties that
	such Dehn extensions are all subgroups of $\SL_2(\CC)$, and that the finitely-generated case
	can be reduced to the finitely-presented case, (Proposition \ref{thruLimit}).
	
	\bigskip\noindent\textbf{Acknowledgement}. 
	The author thanks Ian Agol and Daniel Groves for helpful conversations.

\section{Preliminaries}\label{Sec-Prelim}
	In this section, we recall some notions and results related to the topic of this paper.

	\subsection{The J{\o}rgensen-Thurston Theorem}\label{Subsec-JT}
		For convenience of our discussion, we adopt the following notations
		and state the results in terms of Kleinian groups. Cf.~\cite{Th} for
		facts mentioned in this subsection. 
		
		Let $\Gamma$ be a torsion-free Kleinian group of finite covolume, namely, which has a fundamental domain 
		in $\HH^3$ of finite volume. Then $\Gamma$ has at most finitely many of conjugacy-classes parabolic
		subgroups, represented by subgroups $P^1,\cdots,P^q$, ($q\geq0$), each isomorphic to a rank-$2$
		free abelian group. We often call these chosen subgroups 
		\emph{cusp representatives} of $\Gamma$. By a \emph{slope-tuple} $\zeta=(\zeta^1,\cdots,\zeta^q)$
		with respect to $(\Gamma;P^1,\cdots,P^q)$ (or ambiguously, w.r.t.~$\Gamma$), we shall mean that
		for each $1\leq j\leq q$, the $j$-th component $\zeta^j$ 
		is either trivial or a primitive element in $P^j$. 
		For any slope-tuple $\zeta$ of $\Gamma$, we denote the \emph{Dehn filling}
		of $\Gamma$ along $\zeta$ as:
			$$\Gamma_\zeta=\Gamma\,/\,\overline{\langle\zeta^1,\cdots,\zeta^q\rangle},$$
		i.e.~$\Gamma$ quotienting out the normal closure of the $\zeta^j$'s,
		and often denote the quotient epimorphism as:
			$$\iota_\zeta:\Gamma\twoheadrightarrow\Gamma_\zeta.$$ 
		
		By Thurston's
		Hyperbolic Dehn-Filling Theorem (\cite[Theorem 5.8.2]{Th}), 
		for generic slope-tuples $\zeta$ (avoiding finitely
		many primitive choices for each component), $\Gamma_\zeta$ is isomorphic
		to a torision-free Kleinian group of finite covolume, indeed, strictly less
		than that of $\Gamma$ if $\zeta$ is nonempty.
		These $\Gamma_\zeta$'s are usually called \emph{hyperbolic Dehn fillings} of $\Gamma$. 
		In fact, one may choose faithful Kleinian representations
		$\rho_\zeta:\Gamma_\zeta\to\mobgp$, so that for any sequence
		$\seq{\zeta_n}$ of slope-tuples, there is a subsequence for which the induced
		representations $\seq{\rho_{\zeta_n}\circ\iota_{\zeta_n}}$ of $\Gamma$ strongly converges.
		Moreover, if for each $j$-th component, $\seq{\zeta^j_n}$ has no bounded subsequence
		of primitive elements in $P^j$, then $\seq{\rho_{\zeta_n}\circ\iota_{\zeta_n}}$ strongly converges
		to the inclusion $\Gamma\subset\mobgp$.
		
		The following rephrases \cite[Theorem 5.12.1]{Th}:
		
		\begin{theorem}[{J{\o}rgensen-Thurston}]\label{JT}
			For any $V>0$, there exist finitely many
			torsion-free Kleinian groups $\Gamma_1,\cdots,\Gamma_k$ of covolume at most $V$,
			(together with chosen
			cusp representatives), such that any torsion-free Kleinian group of covolume at most $V$ is isomorphic
			to the hyperbolic Dehn filling $(\Gamma_i)_\zeta$ of some $\Gamma_i$ along some slope-tuple
			$\zeta$ of $\Gamma_i$.
		\end{theorem}
		
		\begin{remark} 
			It is also implied from the proof that hyperbolic Dehn fillings decreases the covolume.
		\end{remark}
		
		\subsection{Dehn extensions}\label{Subsec-DehnExt}
			
			For aspherical orientable compact $3$-manifolds, Dehn extensions
			have been investigated in \cite{AL} from a topological perspective.
			This kind of construction was introduced earlier to define 
			knot invariants, known as generalized knot groups, (cf.~\cite{Ke,LN,Wa}).
			In this note, we rephrase Dehn extensions on the group level
			in terms of amalgamations. We shall restrict ourselves 
			to torsion-free Kleinian groups of finite covolume, for simplicity,
			but we need a multi-cusp version, allowing slope-tuples and denominator-tuples.
			
			Let $(\Gamma;P^1,\cdots,P^q)$ be a torsion-free Kleinian group of finite covolume, 
			together with chosen cusp representatives. As each $P^j<\Gamma$ is isomorphic to a rank-$2$
			free abelian group, we may also naturally identify $P^j$ 
			as the integral lattice of the $\QQ$-vector space $P^j\otimes_\ZZ\QQ$.
			
\begin{figure}[htb]
				\centering
				\psfrag{a}{\scriptsize{$P^1$}}	\psfrag{A}{\scriptsize{$P^1+\ZZ\,\frac{\zeta^1}{m^1}$}}
				\psfrag{b}{\scriptsize{$P^2$}}	\psfrag{B}{\scriptsize{$P^2+\ZZ\,\frac{\zeta^2}{m^2}$}}
				\psfrag{c}{\scriptsize{$P^q$}}	\psfrag{C}{\scriptsize{$P^q+\ZZ\,\frac{\zeta^q}{m^q}$}}
				\psfrag{G}{\scriptsize{$\Gamma$}}
				\includegraphics[scale=.7]{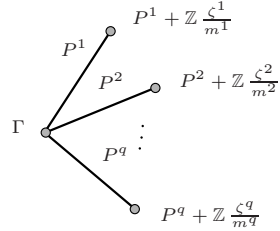}
				\caption{The Dehn extension $\Gamma^{e(\zeta,m)}$.}\label{figTupleDehnExt}
			\end{figure}

			\begin{definition}\label{DehnExt} 
				With the notations above, let $m=(m^1,\cdots,m^q)$ be a $q$-tuple of
				positive integers, called a \emph{denominator-tuple}, and let $\zeta=(\zeta^1,\cdots,\zeta^q)$
				be a slope-tuple, where each $\zeta^j\in P^j$ is either trivial or primitive. 
				The \emph{Dehn extension} of $\Gamma$ along the slope-tuple $\zeta$ with the denominator-tuple
				$m$, denoted by $\Gamma^{e(\zeta,m)}$, is defined as the amalgamation of $\Gamma$ with all the $(P^j+\ZZ\,\frac{\zeta^j}{m^j})$'s (as
				subgroups of $P^j\otimes_\ZZ\QQ$'s)
				along all the $P^j$'s respectively. In other words,
				it is the fundamental group of the graph of groups shown in Figure \ref{figTupleDehnExt}, which we
				write as: 
					$$\Gamma^{e(\zeta,m)}=
						\Gamma\cdot\left(*_{P^1}\,(P^1+\ZZ\,\frac{\zeta^1}{m^1})\right)\cdots\left(*_{P^q}\,(P^q+\ZZ\,\frac{\zeta^q}{m^q})\right).$$
				We often briefly write $\Gamma^e$ whenever $m$ and $\zeta$ are clear from the context.
				The Dehn extension is said to be \emph{trivial on the $j$-th cusp}, if either $\zeta^j$ is
				trivial or $m^j=1$, and it is \emph{trivial} if so it is on each cusp.
			\end{definition}
			
			\begin{definition}\label{cong-dom}
				Let $\zeta,\zeta'$ be two slope-tuples, we say that they are \emph{congruent} modulo
				a denominator-tuple $m$, if for each $1\leq j\leq q$, $\zeta^j\equiv(\zeta')^j\bmod m^j$.
				We say $\zeta$ \emph{dominates} $\zeta'$ modulo $m$, if for each $1\leq j\leq q$,
				either $\zeta^j\equiv(\zeta')^j\bmod m^j$, or that $\zeta^j=0$.
			\end{definition}
			
			It is clear that $\Gamma^{e(\zeta,m)}$ is naturally isomorphic to $\Gamma^{e(\zeta',m)}$, if $\zeta$ is
			congruent to $\zeta'$ modulo $m$. Thus, for any denominator-tuple $m$, there are only finitely
			many distinct Dehn extensions up to isomorphism.
			
			\begin{definition}\label{extFilling}
				Suppose $\Gamma^e$ is the Dehn extension of $\Gamma$ along the slope-tuple $\zeta$
				with the denominator-tuple $m$, then there is a canonical \emph{extended-filling} epimorphism:
					$$\iota^e:\Gamma^e\twoheadrightarrow \Gamma_\zeta,$$
				quotienting out the normal closure of all the $\frac{\zeta^j}{m^j}$'s. Moreover, if $\zeta'$
				is congruent to (resp.~dominated by) $\zeta$ modulo $m$, there are also \emph{congruent} (resp.~\emph{dominated})
				extended-fillings:
					$$\iota^e_{\zeta'}:\Gamma^e\twoheadrightarrow \Gamma_{\zeta'},$$
				quotienting out the normal closure of all the $\frac{(\zeta')^j}{m^j}$'s.
			\end{definition}

\section{The finitely presented case}\label{Sec-fpCase}
	In this section, we prove the finitely presented case of Theorem \ref{main-JTh}, namely:
	
	\begin{prop}\label{JTh-fp} 
		Suppose $G$ is a finitely presented group, and suppose $V>0$. Then
		the conclusion of Theorem \ref{main-JTh} holds for $G$ and $V$.
	\end{prop}

	We prove Proposition \ref{JTh-fp} in the rest of this section.
	It is a consequence of \cite[Theorem 3.2]{AL} and the classical J{\o}rgensen-Thurston
	Theorem.

	\begin{lemma}\label{thruDehnExt}
		Let $G$ be a finitely presented group, and $\Gamma$ be a torsion-free Kleinian group
		of finite covolume, with cusp representatives $P^1,\cdots, P^q$. Then there are finite
		collections of (primitive) slopes $\mathcal{S}^j\subset P^j$ for $1\leq j\leq q$, 
		and there is some denominator-tuple $m=(m^1,\cdots,m^q)$, such that for any slope-tuple $\zeta=(\zeta^1,\cdots,\zeta^q)$
		where each $\zeta^j\not\in\mathcal{S}^j$, and for any homomorphism
		$\phi:G\to\Gamma_\zeta$, $\phi$ factors through the canonical extended-filling epimorphism
		$\iota^e:\Gamma^{e(\zeta,m)}\twoheadrightarrow \Gamma_\zeta$.
	\end{lemma}
	
	\begin{proof}
		By Thurston's Hyperbolic Dehn Filling Theorem (\cite[Theorem 5.8.2]{Th}), we regard $\Gamma$ and all the hyperbolic Dehn fillings
		$\Gamma_\zeta$'s as subgroups 
		of $\mobgp$, so that the Dehn filling epimorphisms $\iota_\zeta:\Gamma\twoheadrightarrow\Gamma_\zeta$ converges strongly
		for any convergent sequence (possibly to $\infty$) of slope-tuples $\zeta$'s, (cf.~Subsection \ref{Subsec-JT}). 
		Choose a base-point $O$ of the
		$3$-dimensional hyperbolic space $\HH^3$
		and an orthonormal basis of the tangent space to identify $\mobgp$ with $\Isom_+(\HH^3)$.
		The quotient manifolds $N_\zeta=\HH^3\,/\,\Gamma_\zeta$
		are all complete finite-volume hyperbolic manfolds with naturally induced base points,
		forming a sequentially-compact family converging to $N_\star=\HH^3\,/\,\Gamma$ 
		in the pointed Gromov-Hausdorff
		topology.
		Let $\epsilon>0$ be a Margulis number of $\HH^3$. Then there are finite collections 
		$\mathcal{S}^j\subset P^j$ of primitive slopes, such that whenever 
		$\zeta$ is a slope-tuple whose components avoid
		elements of the $\mathcal{S}^j$'s, 
		the $\epsilon$-thin part of $N_\zeta$ is the disjoint union
		of $V^1_\zeta,\cdots V^q_\zeta$, each homeomorphic
		to either a solid torus $S^1\times D^2$ or $T^2\times[0,+\infty)$, 
		corresponding to (possibly null) Dehn fillings
		of the cusps.
		Furthermore, one may pick $\mathcal{S}^j$ sufficiently large
		so that to ensure that
		the meridinal area of 
		every solid torus $V^j_\zeta$ is sufficiently large,
		for instance, greater than $A(\ell(G))$ where $\ell(G)$ is the presentation length
		of $G$ (\cite[Definition 3.1]{AL}) and $A(n)=27^n(9n^2+4n)\pi$.

		Note that the argument of \cite[Theorem 3.2]{AL} only used the local geometry of
		the hyperbolic piece near a deep Margulis tube. We may apply the same argument 
		simutaneously to all the tubes $V^1_\zeta,\cdots,V^q_\zeta$ for each $N_\zeta$
		to see that there is some uniform
		denominator-tuple $m=(m^1,\cdots,m^j)$, such that every map $f:K\to N_\zeta$
		factors through the (topological) Dehn extension $N^{e(\zeta,m)}$ via the extended-filling
		map $N^{e(\zeta,m)}\to N_\zeta$, up to homotopy. Here $K$ is a the presentation
		complex of $G$ with $\pi_1(K)\cong G$ (cf.~\cite[Subsection 3.1]{AL}), 
		and $N^{e(\zeta,m)}$ is a natural topological space
		with $\pi_1(N^{e(\zeta,m)})\cong\Gamma^{e(\zeta,m)}$, obtained by attaching `ridge pieces'
		to the cusp tori of $N_\star$ (cf.~\cite[Subsection 2.2]{AL}). In fact,
		one may require each $m^j$ satisfy $0<m^j\leq T(\ell(G))$, where $T(n)=2\cdot 3^n$.
		In language of groups, this is the conclusion of Lemma \ref{thruDehnExt}.
	\end{proof}

	The lemma below may be regarded as a homomorphism analogue to
	Thurston's Hyperbolic Dehn Filling Theorem.

	\begin{lemma}\label{inductionJTh}
		Suppose $G$ is finitely presented, and $\Gamma$
		is a torsion-free Kleinian group of finite covolume. Then there are
		finitely many groups $R_1,\cdots,R_l$, each being a Dehn extension of a hyperbolic
		Dehn filling of $\Gamma$, such that for any hyperbolic Dehn filling
		$\Gamma_\zeta$ of $\Gamma$, there is an $R_i$ an extended filling epimorphism
		$\iota_{i,\zeta}:R_i\twoheadrightarrow\Gamma_\zeta$, so that any homomorphism $\phi:G\to\Gamma_\zeta$
		factors through $\iota_{i,\zeta}$. 
	\end{lemma}
	
	\begin{proof} 
		We make an induction argument on the number of cusps of $\Gamma$.
		When $\Gamma$ has no cusp, there is no further Dehn filling, 
		so there is nothing to prove. Suppose that we have proved the statement 
		when $\Gamma$ has $q-1$ cusps, where $q>0$. 
		
		When $\Gamma$ has $q$ cusps, represented by $P^1,
		\cdots,P^q$, by Lemma \ref{thruDehnExt}, there
		are finite collections of primitive slopes, 
		$\mathcal{S}^j\subset P^j$ for $1\leq j\leq q$,
		and there is a denominator-tuple $m=(m^1,\cdots,m^{q})$,
		such that for any slope-tuple $\zeta$ whose components
		avoid elements of $\mathcal{S}^1,\cdots,\mathcal{S}^{q}$,
		any $\phi:G\to \Gamma_\zeta$ factors through the extended filling
		epimorphism $\iota^e_{\zeta}:\Gamma^{e(\zeta,m)}\twoheadrightarrow \Gamma_\zeta$.
		Note that since $m$ is fixed, there only finitely many isomorphism classes
		of $\Gamma^{e(\zeta,m)}$'s, and two extended fillings $\iota^e_\zeta$ and
		$\iota^e_{\zeta'}$ are congruent when $\zeta$ is congruent to $\zeta'$ modulo
		$m$, (cf.~Subsection \ref{Subsec-DehnExt}). We denote 
		these isomorphically distinct Dehn extensions as $R_1,\cdots,R_{l_0}$.
		
		If $\zeta$ does not satisfy the condition above, it is dominated
		by a slope-tuple $\zeta'$ of the form $(0,\cdots,0,\zeta^j,0,\cdots,0)$ 
		where $\zeta^j\in\mathcal{S}^j$, (cf.~Definition \ref{cong-dom}).
		We enumerate these $\zeta'$'s as $\zeta'_1,\cdots,\zeta'_s$.
		For each $1\leq t\leq s$, since $\Gamma_{\zeta'_t}$ has only $q-1$ cusps, 
		we apply the induction assumption so there are Dehn extensions
		of hyperbolic Dehn fillings of $\Gamma_{\zeta'_t}$, 
		say $R^t_1,\cdots,R^t_{l_t}$, such that if $\zeta$ is dominated 
		by $\zeta'_t$, any $\phi:G\to\Gamma_{\zeta'_t}$ factors through
		the extended filling epimorphism from some $R^t_i$ to $\Gamma_{\zeta'_t}$.
		
		Now we take $R_1,\cdots, R_l$ as $R_1,\cdots,R_{l_0}$ together with
		all the $R^t_i$'s, where $1\leq i\leq l_t$ and $1\leq t\leq s$.
		Then the conclusion of the statement holds for $\Gamma$. This completes the induction.
	\end{proof}
	
	\begin{proof}[Proof of Proposition \ref{JTh-fp}]
		Let $G$ be a finitely presented, and $V>0$ as assumed. By the J{\o}rgensen-Thurston
		Theorem (cf.~Theorem \ref{JT}), there are finitely many torsion-free 
		Kleinian groups $\Gamma_1,\cdots,\Gamma_s$ of covolume at most $V$,
		such that any torsion-free Kleinian group $H$ of covolume at most $V$ 
		is isomorphic to a Dehn filling $(\Gamma_i)_\zeta$ of some $\Gamma_i$
		along some slope-tuple $\zeta$ (with respect to chosen cusp representatives). 
		For each $\Gamma_i$, we apply Lemma \ref{inductionJTh} to obtain 
		finitely many groups $R_{i,1},\cdots,R_{i,k_i}$, each being a Dehn extension
		of a hyperbolic Dehn filling of $\Gamma_i$. Hence all the $R_{i,t}$'s (for all
		$i$'s and $t$'s) are
		Dehn extensions of torsion-free Kleinian groups of covolume at most $V$.
		We rewrite them as $R_1,\cdots,R_l$, then these are the groups as claimed
		by Proposition \ref{JTh-fp}.
	\end{proof}

\section{The finitely generated case}

	In this section, we deduce Theorem \ref{main-JTh} from the finitely presented
	case (Proposition \ref{JTh-fp}). We show that Dehn extensions of 
	finite-covolume torsion-free Kleinian groups are stable limits of
	congruent extended fillings. With some trick of representation varieties,
	this will imply Theorem \ref{main-JTh} using the finitely presented case.

	\subsection{Stable limits}\label{Subsec-stableLim}
		The notion of stable limits is related to the study of limit groups
		in geometric group theory (cf.~\cite{Se,BF}). Our definition here is slightly more general,
		allowing variations on the target groups.
		
		\begin{definition}\label{stableSeq} 
			Let $G$ be a finitely generated group, and $\seq{H_n}$ be a sequence of groups. Suppose $\seq{\phi_n:G\to
			H_n}$ is a sequence of homomorphisms, satisfying that for any $g\in G$,
			either $\phi_n(g)$ is trivial for all but finitely many $n$'s, or that $\phi_n(g)$ is nontrivial for all but finitely
			many $n$'s, then $\seq{\phi_n}$ is said to be a \emph{stable sequence} of homomorphisms.
		\end{definition}
		
		\begin{definition}\label{stableLim} 
			Suppose $\seq{\phi_n:G\to H_n}$ is a stable sequence of homomorphisms, we define the \emph{stable kernel}
			of $\seq{\phi_n}$ to be the (possibly trivial) normal subgroup $K_\infty$ of $G$, consisting of all the elements $g\in G$ each of which is
			trivial under all but finitely many $\phi_n$'s. The quotient epimorphism $\pi:G\twoheadrightarrow G\,/\,K_\infty$,
			(or ambiguously, the quotient group $L_\infty=G\,/\,K_\infty$), is said to be the \emph{stable limit} of the
			stable sequence $\seq{\phi_n}$.
		\end{definition}

	\subsection{Dehn extensions as stable limits}\label{Subsec-extStable}
		In this section, we prove Dehn extensions of torsion-free Kleinian groups
		of finite covolume can be obtained as a stable limit of congruent
		extended fillings. 
		
		Let $\Gamma$ be a torsion-free Kleinian group of finite covolume,
		with cusp representatives $P^1,\cdots,P^q$. We say a slope-tuple
		$\zeta=(\zeta^1,\cdots,\zeta^q)$ is \emph{all-primitive} if every component
		$\zeta^j$ is primitive in $P^j$. An infinite sequence  $\seq{\zeta_n}$
		of all-primitive slope-tuples is said to \emph{converge to all-cusps}, 
		if each component sequence $\seq{\zeta^j_n}$ has no bounded infinite subsequence. 
		
		\begin{prop}\label{extStable} 
			Let $\Gamma$ be a torsion-free Kleinian group of finite covolume,
			with cusp representatives $P^1,\cdots,P^q$. Let $\Gamma^{e(\zeta,m)}$
			be the Dehn extension of $\Gamma$ along some slope-tuple $\zeta$
			with the denominator-tuple $m$, and let $\seq{\zeta_n}$ be a sequence
			of all-primitive slope-tuples, each congruent to $\zeta$ modulo $m$.
			Suppose $\seq{\zeta_n}$ converges to all-cusps. Then $\Gamma^{e(\zeta,m)}$ is the stable
			limit of the stable sequence of congruent extended-filling epimorphisms:
				$$\iota^e_{\zeta_n}:\Gamma^{e(\zeta,m)}\twoheadrightarrow\Gamma_{\zeta_n}.$$ 
		\end{prop}
		
		Proposition \ref{extStable} is a special case of some more general result for Dehn fillings of relatively
		hyperbolic groups. In fact, by a combination theorem for amalgamations of relatively hyperbolic groups
		(cf.~\cite[Theorem 0.1]{Da}, or more generally, \cite[Theorem 1.3]{Os-combination}),
		and by the fact that finitely-generated torsion-free Kleinian groups are hyperbolic relative to
		the cusp subgroups (\cite[Theorem 4.11]{Fa}), Dehn extensions of finite-covolume
		torsion-free Kleinian groups are all hyperbolic relative to the extended cusp subgroups,
		namely, the	$(P^j+\ZZ\,\frac{\zeta^j}{m^j})$'s. Then Proposition \ref{extStable}
		follows immediately from a hyperbolic Dehn filling theorem for relatively
		hyperbolic groups due to Daniel Groves and Jason Manning, (\cite[Corollary 9.7]{GM}, 
		cf.~\cite[Theorem 1.1]{Os-filling}).
		
		For the sake of self-containedness, however, we provide a more direct proof 
		of Proposition \ref{extStable} in the rest of this subsection. 
		Our argument is inspired
		from a paper of Matthew White (\cite{Wh-inj}). The idea is to study
		the asymptotic geometry of any element in the stable kernel $K_\infty$,
		as a word with respect to the defining graph-of-groups decomposition
		of $H^{e(\zeta,m)}$, showing that if it were nontrivial, it would be eventually `fat' 
		in $\HH^3$ and therefore impossible to exist. 
		
		For convenience of notations,
		we shall write the action of $\mobgp$ on $\HH^3$ on the \emph{right}, namely, $x.g$ is the
		action of $g\in\mobgp$ on $x\in\HH^3$; and in a topological space $X$, we write the path-multiplications
		as usual, namely, for a path $\gamma$ from $x$ to $x'$, and a path $\gamma'$ from $x'$ to $x''$,
		$\gamma\gamma'$ means the joined path from $x$ to $x''$.
		
		Suppose the slope-tuple $\zeta=(\zeta^1,\cdots,\zeta^j)$,
		and the denominator-tuple $m=(m^1,\cdots,m^q)$. We write:
			$$\tilde{P}^j=P^j+\ZZ\,\frac{\zeta^j}{m^j},$$
		for $1\leq j\leq q$. Then $\Gamma^e=\Gamma^{e(\zeta,m)}$ is the amalgamations:
			$$\Gamma^e=\Gamma\cdot(*_{P^1}\tilde{P}^1)\cdots(*_{P^q}\tilde{P}^q),$$
		as in Definition \ref{DehnExt}. By conjugations of $\mobgp$,
		we may assume the Dehn filling epimorphisms $\seq{\iota_{\zeta_n}:\Gamma\twoheadrightarrow\Gamma_{\zeta_n}}$
		strongly converges to $\Gamma\subset\mobgp$. We write:
			$$N_\star=\HH^3\,/\,\Gamma, \textrm{ and } N_{\zeta_n}=\HH^3\,/\,\Gamma_{\zeta_n},$$
		with induced base-points $\nu_\star$ and $\nu_n$, respectively.
		Thus $\seq{(N_{\zeta_n},\nu_n)}$ converges to $(N_\star,\nu_\star)$ in the pointed
		Gromov-Hausdorff sense.
		We may also assume the base-points lies in the $\epsilon$-thick parts, where $\epsilon>0$
		is a Margulis number of $\HH^3$. We may suppose
		$\epsilon$ is so small that the $\epsilon$-thin part of $N_\star$ is the disjoint union of
		horocusps of $N_\star$, denoted as:
			$$V_\star=V^1_\star\sqcup\cdots\sqcup V^q_\star.$$ 
		We pick a point $\upsilon^j_\star$ on $\partial V^j_\star$ and let
		$\gamma^j_\star$ be the directed geodesic segment from $\nu_\star$ to $\upsilon^j_\star$
		via which $\pi_1(V^j_\star,\upsilon^j_\star)$ is identified as the subgroup $P^j$ in
		$\pi_1(N_\star,\nu_\star)\cong \Gamma$.
		
		Hence for sufficiently large $n$'s, the $\epsilon$-thin part of $N_{\zeta_n}$ is correspondingly
		a disjoint union of solid tori (by the assumption that $\zeta_n$ is all-primitive), 
		denoted as:
			$$V_n=V^1_n\sqcup\cdots\sqcup V^q_n.$$
		Similarly, there are points $\upsilon^j_n\in\partial V^j_n$, and directed geodesic
		segments $\gamma^j_n$ from $\nu_n$ to $\upsilon^j_n$ for each $1\leq j\leq q$, 
		converging to $\upsilon^j_\star$ and
		$\gamma^j_\star$ respectively, such that $\pi_1(V^j_n,\upsilon^j_n)$ is identified to
		the image of $P^j$ in $\pi_1(N_n,\nu_n)\cong \Gamma_{\zeta_n}$ under the extended filling
		epimorphism $\iota^e_{\zeta_n}$. 
		
		For technical reasons, we also make the
		following convention on the choice of $\epsilon$ ensuring the isolation of $\epsilon$-thin parts:

		\begin{convention}\label{isolated}
			The Margulis number $\epsilon>0$ is so small that the distance between two
			distinct components of the preimage of $V_\star$ in the universal covering $\HH^3$ is at least
			$1+5\delta$, where $\delta={\rm arccosh}(\frac{2}{\sqrt{3}})$ is the upperbound of the inradius of a hyperbolic triangle,
			(realized by the hyperbolic ideal triangle).
		\end{convention}
		
		For simplicity, in this subsection, we shall rewrite the congruent extended fillings $\iota^e_{\zeta_n}$'s as:
			$$\phi_n:\Gamma^e\twoheadrightarrow \Gamma_{\zeta_n}.$$
		Remember we have denoted the stable kernel of $\seq{\phi_n}$ as $K_\infty\leq \Gamma^e$ (Definition
		\ref{stableLim}), and our goal is to show that $K_\infty$ is trivial. For this purpose,
		we must study the hyperbolic geometry of words in $\Gamma^e$ that are killed
		under $\phi_n$. 
		
		For any two elements $g,g'\in H^e$, we write:
			$$g\sim g',$$
		if they are conjugate, i.e.~$g=hg'h^{-1}$ for some $h\in \Gamma^e$. 
		For any $g\in \Gamma^e$, one may write $g$ in a \emph{cyclically reduced normal form} up to conjugacy, 
		with respect to the graph-of-groups decomposition of $\Gamma^e$, namely,
			$$g\sim a_1b_1\cdots a_sb_s,$$
		where the factors $a_i\in \Gamma$, and 
		$b_i\in \tilde{P}^j$ for some $1\leq j\leq q$, and where $s\geq0$ is as small as possible. 
		It is clear that if $a_1b_1\cdots a_sb_s$ is a cyclically reduced normal 
		form, then each $b_i\in \tilde{P}^j\setminus P^j$ for some $1\leq j\leq q$, unless $s=1$ and
		$g$ lies in $\Gamma$ up to conjugacy. Note that the $a_i,b_i$ are only determined up to cyclic permutation
		and elements of $P^j$'s, but in the following we shall usually fix a choice of the cyclically
		reduced normal form.
		
		We start by the following easy observation:
		
		\begin{lemma}\label{nonAmal} 
			$\seq{\phi_n}$ is a stable sequence restricted to either $\Gamma$ or any $\tilde{P}^j$,
			namely, if $g\in \Gamma^e$ is an element of either $\Gamma$ or some $\tilde{P}^j$, then
			$g\in K_\infty$ if and only if $g$ is trivial.
		\end{lemma}
		
		\begin{proof}
			Note that for any $g\in \Gamma$, $\phi_n(g)$ converges to $g$ in $\mobgp$ as $n\to+\infty$.
			It follows from the discreteness that $\seq{\phi_n|_\Gamma}$ stably converges to $\Gamma$
			by some well-known arguments in Kleinian groups.
			
			To see that $\seq{\phi_n|_{\tilde{P}^j}}$ is stable, note $\seq{\phi_n|_{P^j}}$ is already stable as above. 
			For any $g\in \tilde{P}^j$, the image of the $m^j$-th power 
			$\phi_n(g^{m^j})\in \phi_n(P^j)$ is either trivial or nontrivial 
			consistently for all but finitely many $n$'s. 
			As $\phi_n(\tilde{P}^j)$ is a subgroup of $\Gamma_{\zeta^j}$, which is torsion-free, we see
			that $\phi_n(g)$ is consistently trivial or nontrivial 
			for all but finitely many $n$'s as well. This completes the proof.
		\end{proof}
		
		Lemma \ref{nonAmal} allows us to focus on the essential case when $g\in \Gamma^e$ is nontrivially amalgamated.
		
		Now suppose $g\in \Gamma^e$ has been written in a cyclically reduced normal form 
		$a_1b_1\cdots a_sb_s$, where $s\geq1$. 
		In particular, no $b_i$ lies in any $P^j$. To argue by contradiction,
		we suppose futhermore that $g\in K_\infty$.
		
		By Lemma \ref{nonAmal}, given the $g$ as above,
		for any sufficiently large $n$, $\phi_n(b_i)$ are nontrivial for 
		any $1\leq i\leq s$, 
		while $\phi_n(a_1b_1\cdots a_sb_s)=\id$ in $\Gamma_{\zeta_n}$ by the assumption that $g\in K_\infty$. 
		For such an $n$,  in $N_{\zeta_n}$ we may find unique directed loops $\alpha^{(n)}_i$ based at $\nu_n$, 
		and $\beta^{(n)}_i$ based at (the corresponding) $\upsilon^j_n$,
		which are geodesic except at their base-points, 
		representing $\phi_n(a_i)$, $\phi_n(b_i)$ respectively. We often drop the superscript,
		writing $\alpha_i,\beta_i$, if it causes no confusion.
		Note only $\alpha_i$'s could be degenerate. 
		
		Now the word $a_1b_1\cdots a_sb_s$ gives a null-homotopic loop 
		$\alpha_1(\gamma_1\beta_1\bar{\gamma}_1)
		\cdots\alpha_s(\gamma_s\beta_s\bar\gamma_s)$ in $N_{\zeta_n}$ based at $\nu_n$, where $\gamma_i$ is the
		geodesic path $\gamma^j_n$ from $\nu_n$ to $\upsilon^j_n$, according to
		$b_i\in \tilde{P}^j$, and where $\bar\gamma_i$ is its direction-reversal.
		Thus it lifts to $\HH^3$, giving a piecewise-geodesic polygon:
			$$\tilde\alpha_1(\tilde\gamma_1\tilde\beta_1\tilde{\bar\gamma}_1)
				\cdots\tilde\alpha_s(\tilde\gamma_s\tilde\beta_s\tilde{\bar\gamma}_s).$$
		Suppose each $\tilde\beta_i$ is the directed geodesic segment $[x_i,y_i]\subset \HH^3$ (or more precisely,
		$[x^{(n)}_i,y^{(n)}_i]\subset\HH^3$). We consider the piecewise-geodesic polygon $\Pi_n\subset\HH^3$, where:
			$$\Pi_n=[y_s,x_1]\cup[x_1,y_1]\cup\cdots\cup[y_{s-1},x_s]\cup[x_s,y_s].$$
		
		We wish to study the geometry of the polygon $\Pi_n$ for sufficiently large $n$'s.
		
		\begin{lemma}\label{alphaBd} 
			There is a constant $D_1>0$ (depending on $\Gamma$ and the chosen cyclically reduced normal
			form of $g\in \Gamma^e$) such that
			the length of each $[y_{i-1},x_i]$ (where $y_0=y_s$ by convention) is at most $D_1$ for any $n$.
		\end{lemma}
		
		\begin{proof} 
			Since $\seq{(N_{\zeta_n},\nu_n)}$ converges to $(N_\star,\nu_\star)$ in the pointed Gromov-Hausdorff sense, 
			the lengths of $\gamma^j_n$'s and $\alpha_i$'s are uniformly bounded for all $n$'s.
			This implies that the lengths of $[y^{(n)}_{i-1},x^{(n)}_i]$'s are uniformly bounded for all $1\leq i\leq s$
			and all $n$'s.
		\end{proof}
		
		\begin{lemma}\label{betaLong}
			For any constant $L>0$, the length of each $[x_i,y_i]$ becomes greater than
			$L$ when $n$ is sufficiently large (depending on $\Gamma$, $\seq{\zeta_n}$, and 
			the chosen cyclically reduced normal form of $g\in \Gamma^e$).
		\end{lemma}
		
		\begin{proof} 
			To argue by contradiction, supposing that for some $1\leq k\leq s$, there were a subsequence so that the lengths of
			$[x^{(n)}_k,y^{(n)}_k]$'s are uniformly bounded. We may assume $b_k\in \tilde{P}^j$. 
			
			First observe that the translation distances:
				$$d_{\HH^3}(O,\,O.{\phi_n(b_k)}),$$
			are uniformly bounded for all $n$'s, where $O$ is the base-point of $\HH^3$.
			In fact, by our construction, the geodesic $[O,\,O.{\phi_n(b_k)}]$ projects to a loop
			based at $\nu_n\in N_{\zeta_n}$, homotopic to $\gamma^j_n\beta^{(n)}_k\bar\gamma^j_n$ relative
			to $\nu_n$. Thus the observation follows because the length of $\gamma^j_n$ is uniformly
			bounded by the pointed Gromov-Hausdorff convergence, and that the length
			of $\beta^{(n)}_k$ are uniformly bounded by the assumption,
			as it equals that of $[x^{(n)}_k,y^{(n)}_k]$.
			
			Following from this observation, we may passing to a 
			subsequence (still indexed by $n$ for simplicity),
			such that $\phi_n(b_k)$ converges to some
			$b^\star\in \mobgp$. 
			Because $\phi_n(b_k)$ lies in $\phi_n(P^j)=\phi_n(\tilde{P}^j)$, and 
			because $\seq{\phi_n(P^j)}$ converges geometrically $P^j$, 
			we see $b^\star\in P^j$. 
			
			We claim that $b_k=b^\star$ in $\Gamma^e$. To see this, consider $b_k^{-1}b^\star\in \tilde{P}^j$. 
			Since $\phi_n(b_k^{-1}b^\star)\in \Gamma_{\zeta_n}$ converges to the identity in $\mobgp$ by the construction,
			the distances:
				$$d_{\HH^3}(O,\,O.{\phi_n(b_k^{-1}b^\star)}),$$ 
			converges to zero. As
			$O$ projects to $\nu_n\in\partial N_{\zeta_n}$ which 
			lies in the $\epsilon$-thick part, this is possible 
			only if $\phi_n(b_k^{-1}b^\star)$ is trivial for any sufficiently large $n$. By Lemma \ref{nonAmal}, $\seq{\phi_n}$
			is stable restricted to $\tilde{P}^j$, we see $b_k^{-1}b^\star=\id$, so $b_k=b^\star$.
			
			However, we have reached a contradiction as $b_k\not\in P^j$ by our assumption on the cyclically reduced normal
			form, but $b^\star\in P^j$. This completes the proof.
		\end{proof}
		
		\begin{lemma}\label{betaFar} 
			For any sufficiently large $n$ and each $1\leq i\leq s$, the geodesic segments 
			$[x_i,y_i]$ and $[x_{i+1},y_{i+1}]$ (where $x_{s+1}=x_1$ and $y_{s+1}=y_1$ by
			convention) lie in distinct components of the preimage of $V_n\subset N_{\zeta_n}$ in $\HH^3$. 
			Hence for any sufficiently large $n$ (depending on $\Gamma$
			and $\seq{\zeta_n}$), they are at least distance $5\delta$ away from each other,
			according to Covention \ref{isolated}.
		\end{lemma}
		
		\begin{proof} 
			Since $x_i$'s and $y_i$'s project to $\upsilon^j_n$'s in $\partial V_n$, and
			since $V_n\subset N_{\zeta_n}$ is union of the $\epsilon$-thin solid tori in $N_{\zeta_n}$ which are convex, 
			it is clear that any $[x_i,y_i]$ lies entirely in some component of the preimage of $V_n$ in $\HH^3$.
			We must show that the components for $[x_i,y_i]$ and for $[x_{i+1},y_{i+1}]$ are distinct.
			
			To argue by contradiction, suppose both $[x_i,y_i]$ and $[x_{i+1},y_{i+1}]$ 
			lie in a component $\tilde{V}\subset\HH^3$ of the preimage of $V_n$. We may suppose $\tilde{V}$
			projects to the $j$-th component $V^j_n\subset V_n$, so $b_i,b_{i+1}\in \tilde{P}^j$. Note
			$[y_i,x_{i+1}]$ projects to a loop in $N_{\zeta_n}$ based at $\upsilon^j_n\in\partial V^j_n$,
			which is homotopic to $\bar\gamma^j_n\alpha^{(n)}_i\gamma^j_n$ relative to $\upsilon^j_n$. Hence it
			is clear that $x_{i+1}=y_i.\phi_n(a_{i+1})$. This means $\tilde{V}$ meets $\tilde{V}.\phi_n(a_{i+1})$,
			so they coincide for sufficiently large $n$'s by Convention \ref{isolated}. 
			Thus $\phi_n(a_{i+1})\in\phi_n(P^j)$ if $n$ is sufficiently large. Then by a similar argument
			as used in the proof of Lemma \ref{betaLong}, we conclude $a_{i+1}\in P^j$ as well.
			However, this contradicts the assumption on the cyclically reduced normal form,
			as one may replace the sub-word $b_ia_{i+1}b_{i+1}$ by a single element of $\tilde{P}^{j+1}$.
			The `hence' part is straightforward.
		\end{proof}
		
\begin{figure}[htb]
			\centering
			\psfrag{a}{\scriptsize{$y_s$}}	\psfrag{b}{\scriptsize{$x_1$}}	\psfrag{c}{\scriptsize{$z_1$}}	\psfrag{d}{\scriptsize{$y_1$}}
			\psfrag{e}{\scriptsize{$x_2$}}	\psfrag{f}{\scriptsize{$z_2$}}	\psfrag{g}{\scriptsize{$y_2$}}	\psfrag{h}{\scriptsize{$x_s$}}
			\psfrag{i}{\scriptsize{$z_s$}}
			\psfrag{A}{\scriptsize{$\Pi_n$}}
			\psfrag{Z}{\scriptsize{$\hat\Pi_n$}}
			\includegraphics[scale=.7]{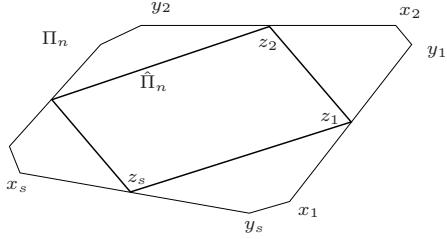}
			\caption{$\Pi_n$ and $\hat\Pi_n$.}\label{figFatPolygon}
		\end{figure}

		Let $\Pi_n$ be the polygon as above, we take $z_i$ to be the midpoints of the corresponding 
		$[x_i,y_i]$ for $1\leq i\leq s$. Let $\hat\Pi_n\subset\HH^3$ be the piecewise-geodesic polygon
		inscribed in $\Pi_n$, where:
			$$\hat\Pi_n=[z_1,z_2]\cup[z_2,z_3]\cup\cdots\cup[z_{s-1},z_s]\cup[z_s,z_1],$$ 
		(Figure \ref{figFatPolygon}). To seek for a contradiction, we show that for any sufficiently large $n$,
		$\hat\Pi_n$ would be too fat to fit in $\HH^3$. In fact, it would have very long sides and 
		would be nearly smooth at the vertices.
		
		We need some terminology at this point. Recall that for any three points $p,q,q'$ in $\HH^3$, the \emph{Gromov
		inner product} between the $q,q'$ with respect to $p$ is defined as:
			$$\langle q,\,q'\rangle_p=\frac12\left(d_{\HH^3}(p,q)+d_{\HH^3}(p,q')-d_{\HH^3}(q,q')\right).$$
		Heuristically, it measures how far from $p$ do the two geodesic segments $[p,q]$ and $[p,q']$
		starts to diverge, so these two segments should `fellow travel' with each other from $p$ until
		reaching this distance.
		
		\begin{lemma}\label{GromovIP} 
			For any $L>0$, the following inqualities hold for sufficiently large $n$'s (depending on
			$N_\star$, ${\zeta_n}$ and the chosen cyclically reduced normal form of 
			$g\in \Gamma^e$) and for each $1\leq i\leq s$: 
				$$\langle z_{i+1},\, y_i\rangle_{z_i}\,>\, \frac{L}{2},$$
			where $z_{s+1}=z_1$ by convention, and:
				$$\langle z_{i-1},\,x_i\rangle_{z_i}\,>\,\frac{L}{2},$$
			where $z_0=z_s$ by convention.
		\end{lemma}
		
		\begin{proof} 
			This follows from some standard coarse geometric estimation. For illustration, we sketch the proof of
			the first inequality. Remember $[x_i,y_i]=\tilde\beta_i$, $[y_i,x_{i+1}]=\tilde\alpha_{i+1}$ and $[x_{i+1},y_{i+1}]
			=\tilde\beta_{i+1}$.
			
			For a given $1\leq i\leq s$, and any $L>0$, let $L_1=2(L+D_1)$, 
			where $D_1>0$ be a constant ensured by Lemma \ref{alphaBd} for sufficiently large $n$'s. 
			By Lemma \ref{betaLong}, for sufficiently large $n$'s, we may assume $d_{\HH^3}(x_i,y_i)> L_1$.
			As $z_i$ is the midpoint of $[x_i,y_i]$,
				$$d_{\HH^3}(z_i,y_i)>\frac{L_1}{2}=L+D_1.$$
			We also assume $n$ is large enough so that Lemma \ref{betaFar} holds. 
			
			Let $v,t,p$ be the points
			on the geodesic segments $[z_i,y_i]$, $[y_i,z_{i+1}]$, $[z_i,z_{i+1}]$ respectively, so that 
			$d_{\HH^3}(z_i,v)=d_{\HH^3}(z_i,p)$, $d_{\HH^3}(y_i,t)=d_{\HH^3}(y_i,v)$, and 
			$d_{\HH^3}(z_{i+1},p)=d_{\HH^3}(z_{i+1},t)$. These are the tangent points of the triangle
			$\Delta(z_i,y_i,z_{i+1})$ with its incircle. Similarly, let $u',v',t'$ be the points on
			the geodesic segments $[y_i,x_{i+1}]$, $[x_{i+1},z_{i+1}]$, $[z_{i+1},y_i]$ respectively,
			so that $d_{\HH^3}(y_i,u')=d_{\HH^3}(y_i,t')$, $d_{\HH^3}(x_{i+1},v')=d_{\HH^3}(x_{i+1},u')$, and 
			$d_{\HH^3}(z_{i+1},t')=d_{\HH^3}(z_{i+1},v')$. See Figure \ref{figHypQuad}.
			
			\begin{figure}[htb]
				\centering
				\psfrag{1}{\scriptsize{$x_i$}}	\psfrag{2}{\scriptsize{$z_i$}}	\psfrag{3}{\scriptsize{$y_i$}}	\psfrag{4}{\scriptsize{$x_{i+1}$}}
				\psfrag{5}{\scriptsize{$z_{i+1}$}}	\psfrag{6}{\scriptsize{$y_{i+1}$}}
				\psfrag{A}{\scriptsize{$u'$}}	\psfrag{b}{\scriptsize{$v$}}	\psfrag{B}{\scriptsize{$v'$}}	
				\psfrag{s}{\scriptsize{$t$}}	\psfrag{S}{\scriptsize{$t'$}}	\psfrag{t}{\scriptsize{$p$}}
				\includegraphics[scale=.7]{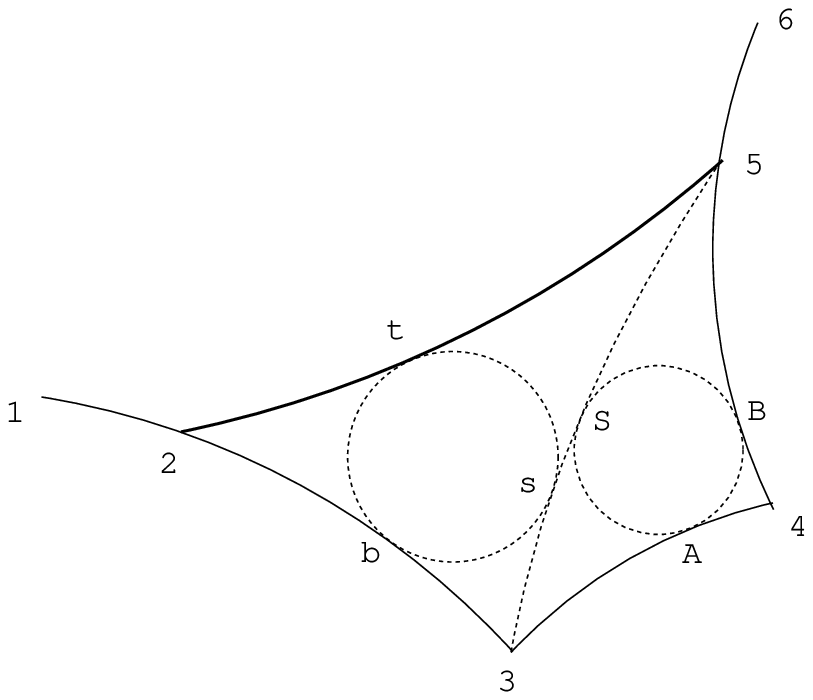}
				\caption{}\label{figHypQuad}
			\end{figure}
			
			We first claim that $t$ lies on $[y_i,t']$. In fact, because:
				$$d_{\HH^3}(v,t)\leq 2\delta,$$
			where $\delta$ is the maximal possible inradius for hyperbolic triangles as in Convention \ref{isolated},
			$[y_i,t]$ lies entirely within the $2\delta$-neighborhood of $[y_i,z_i]$. Similarly, $[t',z_{i+1}]$ lies
			entirely within the $2\delta$-neighborhood of $[x_{i+1},z_{i+1}]$. By the `hence' part of Lemma \ref{betaFar},
			$[y_i,t]$, $[t',z_{i+1}]$ has no overlap, so $t\in [y_i,t']$.
			
			Therefore, 
				$$d_{\HH^3}(y_i,v)=d_{\HH^3}(y_i,t)<d_{\HH^3}(y_i,t')=d_{\HH^3}(y_i,u')<d_{\HH^3}(y_i,x_{i+1})<D_1.$$
			Combining the estimations above:
				$$d_{\HH^3}(z_i,v)=d_{\HH^3}(z_i,y_i)-d_{\HH^3}(y_i,v)>(L+D_1)-D_1=L.$$
				
			As $d_{\HH^3}(z_i,v)$ equals $\langle z_{i+1}\,y_i\rangle_{z_i}$ by definition,
			this proves the first inequality. The second inequality can be proved similarly.
		\end{proof}
		
		Lemma \ref{GromovIP} immediately implies the following:
		
		\begin{lemma}\label{geomC} 
			For any $L>0$ and $0<\theta<\pi$, the side length
			of the polygon $\hat\Pi_n$ is at least $\frac{L}{2}$ and the angles at vertices are at least
			$\pi-\theta$ for sufficiently large $n$'s.
		\end{lemma}
		
		\begin{proof} 
			Apply Lemma \ref{GromovIP}. The side length of $\hat\Pi_n$ are bounded
			below by $\frac{L}{2}$ since $d_{\HH^3}(z_i,z_{i+1})>
			\langle y_i,\, z_{i+1}\rangle_{z_i}$. We may also assume $L$ to be given large enough,
			then the Gromov inner product implies that the angle between $[z_i,y_i]$ and 
			$[z_i,z_{i+1}]$ is less than $\frac{\theta}{2}$, and the same holds for the angle between
			$[z_i,x_i]$ and $[z_i,z_{i-1}]$. The estimation on angles of $\hat\Pi_n$ now follows
			since $z_i$ is the midpoint of the geodesic segment $[x_i,y_i]$.
		\end{proof}
		
		To obtain a contradiction, it suffices to show:
		
		\begin{lemma}\label{tooFat} 
			If $0<\theta<\pi$ and $L>0$ is sufficiently large depending
			only on $\theta$, there is no such polygon $\hat\Pi_n$ in $\HH^3$ satisfying the conclusion of Lemma \ref{GromovIP}.
		\end{lemma}
		
		\begin{proof}
			We consider the diameter $D$ of $\hat\Pi_n\subset\HH^3$. As $\hat\Pi_n$ is a piecewise
			geodesic polygon,
			$D$ is realized the distance between two distinct vertices $z_j$ and $z_k$.
			Let $B$ be the compact ball centered at $z_j$ with radius 
			$D$. Suppose $[z_j,q)$ is the ray in $\HH^3$ that extends $[z_j,z_k]$,
			where $q\in S^2_\infty$ is the intersection with the sphere at infinity.
			Let $w$ be the intersection between the opposite extension of $[z_j,q)$ with $\partial B$. 
			Let $\hat{B}\subset\HH^3$ is the $\HH^3$-closed horoball
			supported at $q$ satisfying $w\in\partial \hat{B}$, (Figure \ref{figWhiteArg}).
			Then $B\subset \hat{B}$,
			and hence $[z_j,z_{j-1}]$, $[z_j,z_{j+1}]$ both lie in $\hat{B}$.
			
			\begin{figure}[htb]
				\centering
				\psfrag{H}{\scriptsize{$\HH^3$}}	\psfrag{S}{\scriptsize{$S^2_\infty$}}	
				\psfrag{T}{\scriptsize{$\hat{B}$}}	\psfrag{B}{\scriptsize{$B$}}
				\psfrag{k}{\scriptsize{$z_k$}}	\psfrag{j}{\scriptsize{$z_j$}}
				\psfrag{1}{\scriptsize{$z_{j-1}$}}	\psfrag{2}{\scriptsize{$z_{j+1}$}}
				\psfrag{q}{\scriptsize{$q$}}	\psfrag{w}{\scriptsize{$w$}}
				\psfrag{C}{\scriptsize{$\hat\Pi_n$}}
				\includegraphics[scale=.7]{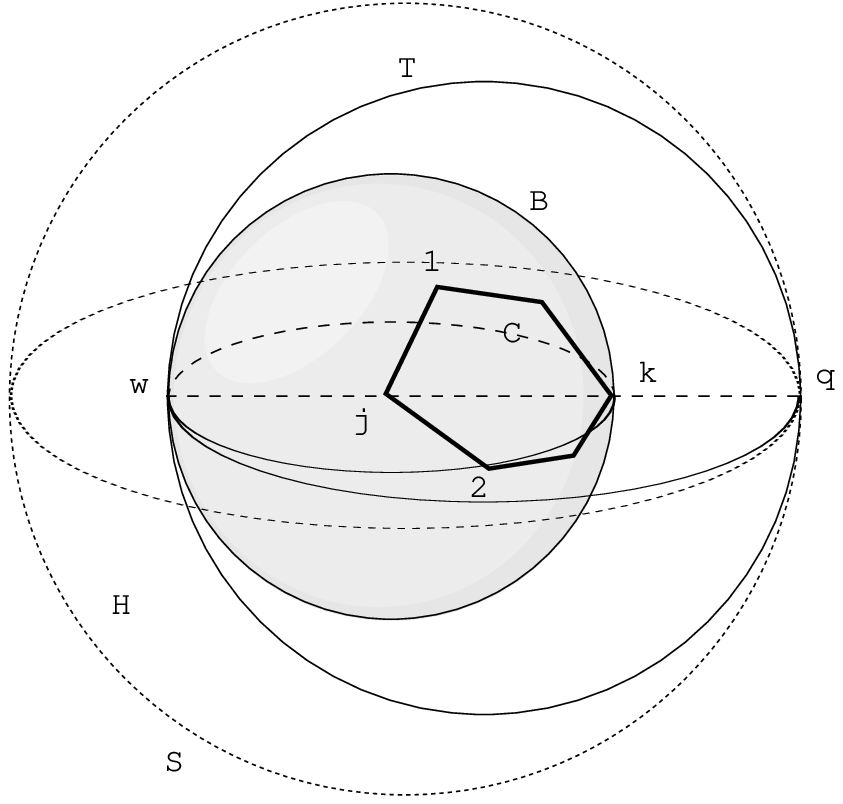}
				\caption{}\label{figWhiteArg}
			\end{figure}
			
			However, since the angle of $\hat\Pi_n$ at $z_j$ is at least $\pi-\theta$, if $L>0$ is sufficiently
			large, at least one of $[z_j,z_{j-1}]$ and $[z_j,z_{j+1}]$ cannot lie entirely in $\hat{B}$.
			Therefore, for such $\theta$ and $L$, a polygon $\Pi_n$ in $\HH^3$ satisfying the
			conclusion of Lemma \ref{GromovIP} would be impossible.
		\end{proof}
		
		To sum up, by Lemma \ref{nonAmal}, the sequence
		$\seq{\phi_n}$ is stable restricted to $\Gamma$ and $\tilde{P}^j$, $1\leq j\leq q$. 
		When $g\in \Gamma^e$ does not lie in $\Gamma$ or
		any $\tilde{P}^j$ up to conjugation, we have shown that $\phi_n(g)$ must be nontrivial for sufficiently
		large $n$'s, since otherwise one would be able to construct a polygon $\hat\Pi_n$ which is forbidden by
		Lemma \ref{tooFat}. Thus $g\not\in K_\infty$ in this case either. 
		We conclude that $\seq{\phi_n:\Gamma^e\twoheadrightarrow
		\Gamma_{\zeta_n}}$ is a stable sequence, and the stable kernel $K_\infty$ is indeed trivial.
		This completes the proof of Proposition \ref{extStable}.

	\subsection{Stable limits of $\SL_2(\CC)$-representations}
		In this subsection, we show that stable limits of a stable sequence
		of $\SL_2(\CC)$-representations are still subgroups of $\SL_2(\CC)$:
				
		\begin{prop}\label{thruLimit}
			Suppose $G$ is a finitely generated
			group, and $L$ is a stable
			limit of a stable sequence of homomorphisms $\seq{\phi_n:G\to H_n}$, where
			each $H_n$ is isomorphic to a subgroup of $\SL_2(\CC)$.
			Then for all but finitely many $n$'s, $\phi_n$
			factors through the stable limit epimorphism
			$\pi:G\twoheadrightarrow L$. 
			Moreover, $L$ is isomorphic a subgroup of $\SL_2(\CC)$.
		\end{prop}
		
		This result follows from some tricks of representation varieties
		which should certainly be familiar to experts of limit groups, (cf.~ 
		\cite{BF,So}). As there seems to be no direct reference in literature,
		we prove Proposition \ref{thruLimit} in the rest of this subsection.
		We refer to
		\cite{CS} for the notion of representation varieties, 
		and \cite{AM} for standard facts in commutative algebra.
		
		\begin{definition}\label{descendant} 
			Let $G$ be a finitely generated group, and $Q$ and $Q'$ are
			quotients of $G$, with
			prescribed quotient maps $\pi,\pi'$. We write $\pi\succ\pi'$, or more often ambiguously:
				$$Q\succ Q',$$ 
			if there is a proper quotient map $\alpha:Q\twoheadrightarrow Q'$ such that
			$\pi'=\alpha\circ\pi$. In this case, $Q'$ is said to be a \emph{descendant} of $Q$.
		\end{definition}
		
		Recall that the representation variety  $R(G)$ 
		of a finitely generated group $G$ in $\SL_2(\CC)$ is the set
		$\Hom(G,\SL_2(\CC))$ with the canonical (Zariski) closed
		affine algebraic structure. 
		More precisely, by picking a generating set $g_1,\cdots,g_r$ of
		$G$, $R(G)$ is identified as a closed affine algebraic subset of
		$\SL_2(\CC)\times\cdots\times\SL_2(\CC)$ ($r$ copies),
		where $\rho\in R(G)$ is identified as $(\phi(g_1),\cdots,\phi(g_r))$. Different
		choices of generators result in natural isomorphisms between algebraic sets.
		
		For any quotient $\pi:G\twoheadrightarrow Q$, there is an induced embedding of variety: 
			$$\pi^*:R(Q)\hookrightarrow R(G),$$ 
		defined by $\pi^*(\rho)=\rho\circ\pi$, for $\rho\in R(Q)$.
		In particular, for every $\rho\in R(G)$,
		the subvariety $R(\rho)\subset R(G)$, is defined as $\rho^*R({\rm Im}(\rho))$. 
		
		The lemma below may be regarded as a representation-variety version 
		of finitely-presented approximation.
		
		\begin{lemma}\label{fpRV}
			For any finitely generated group $G$, there is a 
			finitely presented group $G_0$ of the same rank, (i.e.~the smallest number of generators),
			and a quotient epimorphism $p:G_0\twoheadrightarrow G$ 
			such that $p^*:R(G)\to R(G_0)$ is isomorphic.
		\end{lemma}
		
		\begin{proof}
			Let $\mathcal{P}=(x_1,\cdots,x_r;\,w_1,w_2,\cdots)$ be a (possibly infinite) presentation of $G$,
			where each relator $w_j=w_j(x_1,\cdots,x_r)$ is a word in the generators. Let: 
				$$\xi_i=\left(\begin{array}{cc} a_i&b_i\\c_i&d_i\end{array}\right),$$
			for $1\leq i\leq r$. Then $R(G)\,\cong\,{\rm Spec}^{\rm max}\left(\CC[a_1,\cdots,d_r]\,/\,\mathfrak{a}\right)$,
			where $\mathfrak{a}$ is the ideal generated by all $\det(\xi_i)-1$, and all the entries of $w_j(\xi_1,\cdots,\xi_r)$.
			By Hilbert's Basis Theorem, $\mathfrak{a}$ is finitely generated, say by $u_1,\cdots u_m$, each
			of which is a polynomial in $\det(\xi_i)-1$ and entries of $w_j(\xi_1,\cdots,\xi_r)$. Let $w_{j_1},\cdots
			w_{j_k}$ be the finite set of relators of which at least one entry appears in the polynomials above. 
			Pick $G_0=\langle x_1,\cdots,x_r;\, w_{j_1},\cdots,w_{j_k}\rangle,$
			and the natural quotient epimorphism $p:G_0\twoheadrightarrow G$. Then $R(G_0)$ and $R(G)$ are defined by the same set
			of polynomials in $\CC^{4r}$, and hence $p^*:R(G)\to R(G_0)$ is isomorphic.
		\end{proof}
		
		\begin{lemma}\label{genericRep} 
			Every irreducible component of $R(G)$ contains a generic representation $\rho$, in the sense that
			the embedded image of $R(\rho)$ contains the component. Moreover, any two generic
			representations on the same component have the same kernel in $G$.
		\end{lemma}
		
		\begin{proof} 
			Let $G$ be $r$-generated, so by picking a set of generators,
			we may identify $R(G)$ as embedded in $R(\FF^r)$, where $\FF^r$ is the free group
			of rank $r$. Let $X$ be any irreducible component of $R(G)$. 
			For any $\rho\in X$, by Lemma \ref{fpRV}, there is a finitely presented
			quotient $Q_0$ of $\FF^r$, such that the image
			of $R(Q_0)\hookrightarrow R(\FF^r)$
			coincides with that of $R(\rho)\hookrightarrow R(\FF^r)$. Thus there
			are at most countably many $R(\rho)\hookrightarrow R(\FF^r)$ with distinct
			images. By Baire's Category Theorem, there must be some $\rho\in X$
			such that $\dim R(\rho)=\dim X$, so the embedded image of 
			$R(\rho)$ contains $X$. The
			`moreover' part follows because the representation varieties of two generic
			representations contain each other.
		\end{proof}
		
		For $\rho,\rho'\in R(G)$, we write $\rho\succ\rho'$ meaning that ${\rm Im}(\rho)\succ {\rm Im}(\rho')$.
		
		\begin{lemma}\label{dccRV} 
			Let $G$ be finitely generated, then every proper descending chain of $\SL_2(\CC)$
			representations terminates after finitely many steps.
		\end{lemma}
		
		\begin{proof}
			Suppose $\rho_1\succ\rho_2\succ\cdots$. Note because for $R(\rho_{n+1})\subset R(\rho_n)$, 
			either $R(\rho_{n+1})$ has lower dimension (by `\emph{dimension}' we mean
			the highest dimension among irreducible components) than 
			$R(\rho_n)$ does, or they have the same dimension, but $R(\rho_{n+1})$ has fewer irreducible
			components than $R(\rho_n)$ does. Thus the sequence must terminate at some step.
		\end{proof}
		
		Now we are ready to prove Proposition \ref{thruLimit}.
		
		\begin{proof}[Proof of Proposition \ref{thruLimit}]
			Let $R(G)$ be the representation variety of $G$ in $\SL_2(\CC)$.
			By choosing an identification of each $H_n$ to a subgroup of $\SL_2(\CC)$, we regard each $\phi_n:G\to H_n$
			as a representation $\rho_n:G\to\SL_2(\CC)$. To argue by contradiction, suppose there were an infinite
			subsequence, still indexed by $n$, such that no $\rho_n$ factors through $\pi:G\twoheadrightarrow L$.
			
			Let $\mathcal{Q}$ be the set of representations $\rho\in R(G)$ which descends to infinitely many $\rho_n$'s. 
			The set $\mathcal{Q}$ is non-empty as it contains the generic 
			$\rho\in X$ of an irreducible component $X\subset R(G)$ (by Lemma \ref{genericRep})
			on which there are infinitely many $\rho_n$'s. By Lemma \ref{dccRV}, any proper
			descending sequence of representations terminates, so there is a minimal
			element $\rho^*:G\to \SL_2(\CC)$ in $\mathcal{Q}$ under `$\succeq$'. By our assumption,
			$\rho^*\succeq \rho_{n'}$ for an infinite subsequence $\seq{\rho_{n'}}$. 
			As $\seq{\rho_{n'}}$ is still a stable sequence, clearly 
			${\rm Ker}(\rho^*)\subset K_\infty$, where $K_\infty$
			denotes the stable kernel. Pick an element $g\in K_\infty$, and
			let $Q$ be ${\rm Im}(\rho^*)$ further quotienting out $\rho^*(g)$. 
			Then the embedded image of $R(Q)\hookrightarrow R(\rho^*)$ also contains
			infinitely many $\rho_{n'}$'s. 
			Take a generic representation $\rho^{**}\subset R(G)$ of an irreducible component of 
			the embedded image of $R(Q)\hookrightarrow R(G)$ carrying infinitely
			many $\rho_{n'}$'s, 
			then $\rho^*\succeq\rho^{**}$. It follows that
			$\rho^*=\rho^{**}$ by minimality of $\rho^*$ in $\mathcal{Q}$. 
			Thus $\rho^*(g)=\rho^{**}(g)=\id$. 
			This implies $K_\infty\subset{\rm Ker}(\rho^*)$. 
			Therefore, $K_\infty={\rm Ker}(\rho^*)$, so the subsequence $\seq{\phi_{n'}}$ 
			factors through $\pi:G\twoheadrightarrow L$
			as $\rho^*\succeq\rho_{n'}$ and $L=G\,/\,K_\infty$. This is contrary to our
			assumption, so we have proved the `factorization' part. 
			
			The `moreover' part also follows as ${\rm Im}(\rho_0)\cong L$.
		\end{proof}

	\subsection{Proof of the finitely generated case}
		We now give the proof of Theorem \ref{main-JTh}.
		
		\begin{proof}[Proof of Theorem \ref{main-JTh}]
			Suppose $G$ is a finitely generated group, and $V>0$. By Lemma \ref{fpRV},
			there is a finitely presented group $G_0$ and an epimorphism $p:G_0\twoheadrightarrow G$,
			which induces an isomorphism $p^*:R(G)\to R(G_0)$ on the $\SL_2(\CC)$-varieties.
			By Proposition \ref{JTh-fp}, we may find finitely many Dehn extensions
			of torsion-free Kleinian groups of covolume at most $V$, say $R_1,\cdots,R_k$,
			and for any torsion-free Kleinian group $H$, there is some $R_i$ and an extended-filling
			epimorphism $\iota^e_{i,H}:R_i\twoheadrightarrow H$, such that any homomorphism $\phi:G_0
			\to H$ factors through $\iota^e_{i,H}$. 
			
			We claim these $R_i$'s (and $\iota^e_{i,H}$'s) 
			also satisfy the conclusion of Theorem \ref{main-JTh}.
			A result of William Thurston asserts that finitely-generated torsion-free Kleinian groups
			can be lifted to be embedded in $\SL_2(\CC)$, (cf.~\cite[Proposition 3.1.1]{CS}). 
			Thus, by Propositions \ref{extStable} and \ref{thruLimit}, there are faithful
			representations $\rho_i:R_i\hookrightarrow\SL_2(\CC)$. To see the
			claim, consider any homomorphism $\phi:G\to H$, where $H$ is as above.
			Then $\phi\circ p:G_0\to H$ factors through some $\iota^e_{i,H}$, say $\phi\circ p=\iota^e_{i,H}\circ\psi$,
			where $\psi:G_0\to R_i$. Thus $\rho_i\circ\psi:G_0\to\SL_2(\CC)$ is a representation
			of $G_0$ in $\SL_2(\CC)$. Because $p^*:R(G)\to R(G_0)$ is isomorphic, $\rho_i\circ\psi$
			descends to a representation of $G$ via $p$, or equivalently, 
			$\psi=\bar{\psi}\circ p$ for some $\bar\psi:G\to {\rm Im}(\psi)\hookrightarrow R_i$.
			Hence $\phi\circ p=\iota^e_{i,H}\circ\bar\psi\circ p$, so $\phi=\iota^e_{i,H}\circ\bar\psi$
			as $p$ is surjective. This means $\phi$ factors through $\iota^e_{i,H}$ via $\bar\psi$.
			Then we have proved the claim. This also completes the proof of Theorem \ref{main-JTh}.
		\end{proof}

\section{Conclusions}
	In conclusion, Dehn extensions arise naturally from a limit process
	in the study of homomorphisms to torsion-free
	Kleinian groups with uniformly bounded covolume. From the perspective 
	of geometric group theory, it seems reasonable to expect that there is 
	a similar version of Theorem \ref{main-JTh} describing the
	structure of homomorphisms from a finitely generated group to the Dehn-filling family of a 
	torsion-free relatively hyperbolic group with abelian flats. However,
	as it relies severely on the fact that Kleinian groups are subgroups of $\SL_2(\CC)$, 
	our approach in this note does not 
	apply in general. We also wonder whether it is possible to describe the
	homomorphism structure if we drop the assumption of
	uniformly bounded covolume in Theorem \ref{main-JTh}.

\bibliographystyle{amsalpha}

\end{document}